\documentclass[12pt]{article}

\usepackage{graphicx, psfrag, color, amsmath, amscd, geometry, url}
\input epsf.sty
\epsfverbosetrue 

\usepackage{latexsym}
\usepackage{amsfonts}
\usepackage{amssymb}

\newcommand{\R}{\mathbb{R}}
\newcommand{\Q}{\mathbb{Q}}
\newcommand{\Z}{\mathbb{Z}}
\newcommand{\C}{\mathbb{C}}

\newcommand{\N}{\mathbb{N}}

\newcommand{\rs}{\mbox{$\widehat{\C}$}}

\def\AAA{{\mathcal A}}

\def\CCC{{\mathcal C}}
\def\DDD{{\mathcal D}}

\def\GGG{{\mathcal G}}

\def\JJJ{{\mathcal J}}

\def\MMM{{\mathcal M}}

\def\OOO{{\mathcal O}}
\def\MMM{{\mathcal M}}

\def\YYY{{\mathcal Y}}
\def\NNN{{\mathcal N}}

\def\RRR{{\mathcal R}}

\def\UUU{{\mathcal U}}

\def\sfs{f_\Sigma}

\newtheorem{thm}{Theorem}[section]
\newtheorem{defn}[thm]{Definition}
\newtheorem{remark}[thm]{Remark}

\newtheorem{prop}[thm]{Proposition}

\newtheorem{lemma}[thm]{Lemma}

\newcommand{\qed}{\nopagebreak \begin{flushright}
       \rule{2mm}{2.5mm} \end{flushright}}
\newcommand{\qedspecial}[1]{\nopagebreak \begin{flushright}
       \rule{2mm}{2.5mm}{\bf #1} \end{flushright}}
\newcommand{\bdry}{\partial}                     
\newcommand{\id}{\mbox{\rm id}}                  
\newcommand{\cl}{\overline}                      

\newcommand{\intersect}{\cap}                    
\newcommand{\union}{\cup}                        
\newcommand{\diam}{\mbox{\rm diam}}					         

\newcommand{\lcm}{\mbox{\rm lcm}}      

\newcommand{\mtwo}[4]                            
{\mbox{$\left(\begin{array}{cc}                  
#1 & #2 \\
#3 & #4 
\end{array}
\right)$}}

\newcommand{\dettwo}[4]                          
{\mbox{$\left|\begin{array}{cc}                  
#1 & #2 \\
#3 & #4 
\end{array}
\right|$}}

\newcommand{\pf}{\noindent {\bf Proof: }}

\newcommand{\be}{\begin{enumerate}}
\newcommand{\eb}{\end{enumerate}}
\newcommand{\bi}{\begin{itemize}}
\newcommand{\ib}{\end{itemize}}
\newcommand{\bl}{\begin{list}}
\newcommand{\lb}{\end{list}}

\newcommand{\gap}{\vspace{5pt}}                 

\newcommand{\al}{\alpha} 
\newcommand{\ep}{\varepsilon} 

\newcommand{\wtE}{\widetilde{E}}
\newcommand{\wtU}{{\widetilde{U}}}
\newcommand{\wtV}{\widetilde{V}}

\newcommand{\wtgamma}{\widetilde{\gamma}}

\newcommand{\txi}{\tilde{\xi}}

\renewcommand\tilde{\widetilde}

\def\pullback{\longleftarrow}
\def\univpullback{\hookleftarrow}
\newcommand{\trace}{\mbox{\rm tr}}

\begin{document}

\title{An algebraic characterization of expanding Thurston maps}

\author{Peter Ha\"issinsky and Kevin M. Pilgrim}

\date{\today}

\maketitle

\abstract{Let $f: S^2 \to S^2$ be a postcritically finite branched covering map without periodic branch points.  We give necessary and sufficient algebraic conditions for $f$ to be homotopic, relative to its postcritical set, to an expanding map $g$.}

\tableofcontents

\section{Introduction}

Let $T^2 = \R^2/\Z^2$ denote the usual torus and, using the standard basis, identify $T^2$ with the quotient $H_1(T^2, \R)/H_1(T^2, \Z)$.  Equip $T^2$ with the Euclidean metric.  The following result is well-known:
\begin{thm}
\label{thm:torusversion} 
Suppose $f: T^2 \to T^2$ is a self-covering map of the torus.  
Then $f$ is homotopic through covering maps to an expanding map $g: T^2 \to T^2$ if and only if the spectrum of $f_*: H_1(T^2, \R) \to H_1(T^2, \R)$ lies outside the closed unit disk.  
\end{thm} 
In this setting, one may take $g$ to be the linear map on $T^2=H_1(T^2, \R)/H_1(T^2, \Z)$ induced by $f_*$.  One may further arrange so that the map $g$ is a factor of $f$ via a map $\pi: T^2 \to T^2$ (that is, $\pi \circ f = g \circ \pi$) where $\pi$  is homotopic to the identity,  is monotone, has fibers not disconnecting $T^2$, and, if $f$ is itself expanding, is a homeomorphism.  The proof is now standard: one lifts the identity map $h_0: T^2 \to T^2$ under $f$ and $g$ to obtain $h_1: T^2 \to T^2$ and a homotopy $h_t, t \in [0,1]$ from $h_0$ to $h_1$.  By induction, lifting and concatenating the homotopies, one obtains a family $h_t, t \in [0,\infty)$, such that $g\circ h_{t+1} = h_t \circ f$.  Expansion of $g$ implies the family $h_t$ converges uniformly to the desired semiconjugacy.   Similar and much more general results hold, cf. \cite[Theorem 3]{shub:endomporphisms}. 

In this work, we formulate and prove an analog of this theorem for a certain class of continuous maps $f: S^2 \to S^2$ of the sphere to itself, called {\em Thurston maps}.   

A Thurston map $f: S^2 \to S^2$ is an orientation-preserving branched covering of the two-sphere  of degree $d \geq 2$ for which the necessarily nonempty {\em postcritical set} 
$P_f= P = \union_{n>0}f^{\circ n}(\{\mbox{branch points}\})$ is finite.  These were originally introduced by Thurston \cite{DH1} as combinatorial objects 
with which to classify rational functions as dynamical systems on the Riemann sphere $\rs$.  Two Thurston maps $f, g$ are {\em combinatorially equivalent} 
if there exist orientation-preserving homeomorphisms 
$h_0, h_1: (S^2, P_f) \to (S^2, P_g)$ such that $h_0\circ f = g \circ h_1$ and $h_0, h_1$ are homotopic relative to $P_f$.  
The relation of combinatorial equivalence can be thought of as conjugacy, up to isotopy relative to postcritical sets. 
Recently,  Thurston maps  have been investigated in connection with the study of finite subdivision rules, the classification of quasisymmetry classes 
of metrics on $S^2$, and Cannon's conjecture; see \cite{cfp:fsr}, \cite{bonk:meyer:subdivisions}, \cite{kmp:ph:cxci}.   

A formulation of a variant of Theorem \ref{thm:torusversion} for Thurston maps requires the development of a notion of expansion for non-locally injective maps, 
algebraic/ homotopy-theoretic invariants of a self-map of a simply-connected surface, and the synthetic construction of an expanding model map from the algebraic data. 
\gap

\noindent{\bf Expansion.} Let $f$ be a Thurston map. Since $f$ has branch points, it can never be expanding with respect to a Riemannian metric on $S^2$, nor can it be positively expansive.  We therefore reformulate the notion of expansion in terms of contraction of inverse branches.  Let $\UUU_0$ be a finite cover of $S^2$ by connected open sets, and let $\UUU_n$ denote the covering whose elements are connected components of $f^{-n}(U), U \in \UUU_0$.  We say that $f$ is {\em expanding} if the following condition holds. 

\begin{defn}
\label{defn:expansion} The map $f$ is {\em expanding} if the mesh of the coverings $\UUU_n$ tends to zero as $n \to \infty$.  That is, for any finite open cover $\YYY$ of $S^2$ by open sets, there exists $N$ such that for all $n \geq N$ and all $U \in \UUU_n$, there exists $Y \in \YYY$ with $U \subset Y$.  Equivalently: the diameters of the elements of $\UUU_n$, with respect to some (any) metric, tend to zero as $n \to \infty$.
\end{defn} 

Note that this is a purely topological condition.  It is at present unknown whether every expanding Thurston map is topologically conjugate to a map which is smooth away from a finite set of points.  Perhaps more surprisingly, it is unknown whether every expanding Thurston map is topologically conjugate to a map which preserves some complete length structure on $S^2$.   
\gap

\noindent{\bf Algebraic invariants.}  Nekrashevych \cite{nekrashevych:book:selfsimilar} introduced a collection of algebraic invariants associated to a very wide class 
of dynamical settings, including that of a Thurston map $f$.   These may be unfamiliar to most readers; see \S 3 for details. As in \cite{DH1}, associated to 
$f$ is an orbifold structure $\OOO$ on $S^2$.  In order to have needed finiteness properties, we assume $f$ has no periodic branch points, so that $\OOO$ is compact.  
Once a basepoint is fixed, there is a corresponding orbifold fundamental group $\pi_1(\OOO)$, and a so-called {\em virtual endomorphism} 
$\phi: \pi_1(\OOO) \dasharrow \pi_1(\OOO)$, well-defined up to postcomposition by inner automorphisms.  
Roughly speaking, $\phi$ is the homomorphism induced by the {\em inverse} of $f$.   The property of $\phi$ being {\em contracting} plays 
the role of the spectral condition on the induced map $f_*$ on homology in Theorem \ref{thm:torusversion}. 
\gap

\noindent{\bf Synthetic model map.} When the virtual endomorphism $\phi$ is contracting,  by choosing suitable defining data 
(basis; generators) one may construct a negatively curved {\em self-similarity complex} $\Sigma$.  The construction is similar to that of the boundary 
of a hyperbolic group from a Cayley graph.  The boundary at infinity $\mathcal{J}$ of $\Sigma$ inherits a self-map $\sfs: \JJJ \to \JJJ$.  
The dynamical system $\sfs: \JJJ \to \JJJ$ plays the role of the induced linear map $g$ in Theorem \ref{thm:torusversion}. 
\gap

Our main result is 

\begin{thm}[Characterization of expanding Thurston maps]
\label{thm:characterization} 
Suppose $f$ is a Thurston map without periodic critical points. Then $f$ is homotopic (through Thurston maps relative to its postcritical set) to an expanding map $g$ 
if and only if the virtual endomorphism $\phi: \pi_1(\OOO) \dasharrow \pi_1(\OOO)$ is contracting.  \end{thm}

By applying an observation of Rees \cite[\S\,1.4]{rees:degree2:part1}, we conclude that a Thurston map without periodic branch points is combinatorially equivalent to an expanding map if and only if its virtual endomorphism on the orbifold fundamental group is contracting. 

The hypothesis of no periodic critical points cannot be dropped: the rational function $f(z)=z^2$ acting on the Riemann sphere has contracting virtual endomorphism, but is not homotopic to any expanding map: the nonnegative real axis is an invariant arc joining periodic points in its postcritical set, and by Proposition \ref{prop:obstructions_to_expansion} this is an obstruction to expansion. 

If the defining data for $\Sigma$ and an open covering $\UUU_0$ of $S^2$  is chosen in a suitable way, one obtains a semiconjugacy $\pi: S^2 \to \JJJ$ from 
$f: S^2 \to S^2$ to $\sfs: \JJJ \to \JJJ$ playing the role of the map $\pi$ introduced after Theorem \ref{thm:torusversion}.  In the setting of Theorem \ref{thm:characterization}, 
$\JJJ$ is a sphere, and $g$  is a limit of conjugates of $f$ via a path of conjugacies $h_t$ (with now $t \in [0,1])$ which gradually collapse the fibers of $\pi$ to points. 
Formally, $h_t$ is a {\em pseudo-isotopy}.  

A {\em finite subdivision rule} $\RRR$ on the sphere (in the sense of \cite{cfp:fsr})  is, roughly, a Thurston map $f: S^2 \to S^2$ together with a cell structure 
$S_\RRR$ on $S^2$ and a refinement $\RRR(S_\RRR)$ of $S_\RRR$ such that $f: \RRR(S_\RRR) \to S_\RRR$ is cellular and a homeomorphism on each cell.  
By taking preimages, one obtains a sequence of refinements $\RRR^n(S_\RRR)$, $n \in \N$. One says that $\RRR$ has {\em mesh going to zero} 
if the diameters of the cells at level $n$ tend to zero as $n$ tends to infinity.  It has {\em mesh going to zero combinatorially} if for some $n$, (i) every edge of 
$S_\RRR$ is properly subdivided in $\RRR^n(S_\RRR)$, and (ii) given a pair of disjoint edges $e, e'$ of $S_\RRR$ in the boundary of a tile, 
no tile $t$ of $\RRR^n(S_\RRR)$ contains edges in both $e$ and in $e'$.  Thus mesh going to zero is clearly an invariant of topological conjugacy, 
while mesh going to zero combinatorially is a combinatorial property. 


By appealing to recent work of Bonk and Meyer \cite{bonk:meyer:subdivisions}, we obtain 

\begin{thm}[Equivalence of combinatorial expansion]
\label{thm:equivalence}
Suppose $f$ is a Thurston map without periodic branch points. The following conditions  are equivalent. 
\be
\item The map $f$ is  equivalent to an expanding map, $g$. 
\item There exists $m$ such that $f^m$ is equivalent to a map $g$ which is the subdivision map of a subdivision rule $\RRR$ with mesh going to zero. 
\item There exists $m$ such that $f^m$ is equivalent to a map $g$ which is the subdivision map of a subdivision rule $\RRR$ with mesh going to zero combinatorially. 

\item The virtual endomorphism $\phi_f$ is contracting.
\eb
\end{thm}

\pf That (1) implies (2)  follows from \cite[Theorem 1.2]{bonk:meyer:subdivisions}; we remark that they do not require the absence of periodic branch points.  That (2) implies (1) is a straightforward consequence of the definitions. The equivalence of (2) and (3) is \cite[Theorem 3.2]{cfp:fsr}.   Now suppose (3) holds. By \cite[Theorem 1.1]{cfpp:fsr-cve}, the virtual endomorphism $\phi_g$ is contracting. This implies (4).  That (4) implies (1) is Theorem \ref{thm:characterization}. \qed

In Theorem \ref{thm:equivalence}, the topological conjugacy class of the expanding map $g$ is necessarily unique, by e.g. \cite[Theorem 10.4]{bonk:meyer:subdivisions}.  A corollary of Theorem  \ref{thm:equivalence} is thus that the dynamics of an expanding Thurston map is determined by either of two pieces of combinatorial data: a subdivision rule with mesh going to zero combinatorially, or by an expanding virtual endomorphism.   If a Thurston map $f$ is already presented in terms of a subdivision rule $\RRR$, there is a finite, easily checkable, necessary and sufficient condition for $\RRR$ to have mesh going to zero combinatorially \cite[Theorem 6.1]{cfpp:fsr-cve}.   If $f$ is not so presented, one can try to verify the contraction of $\phi_f$ instead.  Bartholdi has pointed out that this can be done algorithmically, e.g. by combining a normal forms algorithm for writing elements of the orbifold fundamental group with the FR package \cite{bartholdi:FR} in the computational algebra program GAP.

In the paragraphs below, we outline in detail the reduction of the proof of Theorem \ref{thm:characterization} to a few key steps.
\gap

\noindent{\bf Necessity.}  
The implications (1) $\implies$ (2) $\implies$ (3) $\implies$ (4) have already been established in the above proof of Theorem \ref{thm:equivalence}.  In \S 7, however, we provide a new, direct proof.  The idea is that
the group $\pi_1(\OOO)$ acts on the universal cover $p:\tilde{\OOO}\to\OOO$ and that we may find a lift $\tilde{f}_{-}$
of ``$f^{-1}$'' i.e, a map satisfying $f\circ p\circ \tilde{f}_{-}=p$. The action of the virtual endomorphism on $\pi_1(\OOO)$
can be interpreted in terms of the action of $\tilde{f}_{-}$ on $\tilde{\OOO}$. 
Assuming that $f$ is expanding yields a metric on $\OOO$
which can be lifted to $\tilde{\OOO}$ in such a way that (a) $\pi_1(\OOO)$ acts geometrically on $\tilde{\OOO}$
implying that $\pi_1(\OOO)$ is quasi-isometric to $\tilde{\OOO}$ and (b) the action of $\tilde{f}_{-}$ is contracting.
Bringing together these two facts leads to the contraction of the virtual endomorphism.

\gap

\noindent{\bf Sufficiency.} Suppose the virtual endomorphism $\phi$ is contracting and $f$ has no periodic branch points.  Necessarily  $\#P \geq 3$.   For technical reasons, it is convenient later to assume $\#P > 3$. When $\#P=3$, 
$f$ is equivalent to a rational function which is chaotic on the whole sphere, by Thurston's characterization \cite{DH1}, and such a map is expanding.  

We first construct the semiconjugacy $\pi: S^2 \to \JJJ$ from $f: S^2 \to S^2$ to the model map $f_\Sigma: \JJJ \to \JJJ$; recall that $\JJJ$ arises as the boundary of the hyperbolic selfsimilarity complex, $\Sigma$.  Let $\UUU_0$ be an open cover of $S^2$ by Jordan domains.  By taking inverse images, there is another associated infinite one-complex $\Gamma$; see \cite{kmp:ph:cxci} and \S 5.    Using the absence of periodic branch points, we show (Theorem \ref{thm:coincidence}) that for  suitable choices of the covering $\UUU_0$ defining $\Gamma$ and the  defining data for $\Sigma$, there is a natural quasi-isometry $\Phi: \Sigma \to \Gamma$.  Thus $\Gamma$ is also hyperbolic.   Extending $\Phi$ to the boundary by a map of the same name, the composition $\pi:=\Phi^{-1} \circ \pi_\Gamma$ is the desired semiconjugacy.

Next, we show that $\JJJ$ is homeomorphic to $S^2$, and that the model map $f_\Sigma$ is a Thurston map.  A key role is played by the fact, established in \cite{kmp:ph:cxcii}, that $f_\Sigma: \JJJ \to \JJJ$ is a branched covering that, as a topological dynamical system, is expanding and is of so-called {\em finite type}.   In particular, the local degrees of iterates are uniformly bounded.   We first develop (\S 4) some dynamical consequences of the hypothesis  that $\phi$ is contracting, e.g. that there are no so-called Levy cycles. Next, after showing $\pi$ is monotone (Lemma \ref{lemma:monotone}), it suffices, by a classical theorem of Moore  \cite[Theorem 25.1]{daverman:decompositions}, to show  that the fibers of $\pi$ do not disconnect $S^2$.   To see this, we use a case-by-case analysis of the possibilities for how such fibers separate and intersect $P$ and the absence of Levy cycles and other related obstructions to expansion.  We conclude that $\JJJ$ is a sphere, that $f_\Sigma: \JJJ \to \JJJ$ is an expanding Thurston map, and that $\pi|_P$ is injective. 

It remains to construct an isotopy from $f$ to an expanding map $g$. By \cite[Theorems 13.4, 25.1]{daverman:decompositions}, the decomposition $\GGG$ of $S^2$ by the fibers of the semiconjugacy $\pi: S^2 \to \JJJ$ has the property of being {\em strongly shrinkable}: there is a one-parameter family of continuous maps $h_t: S^2 \to S^2, t \in [0,1]$, called a {\em pseudoisotopy}, such that $h_0 = \id_{S^2}$, $h_t$ is a homeomorphism for $t \in [0,1)$, $h_t|_{P_f} = \id_{P_f}$ for all $t$, and the fibers of $h_1$ and those of $\pi$ coincide.  The induced homeomorphism $S^2 \to \JJJ$ conjugates $\sfs$ to an expanding Thurston map $g: S^2 \to S^2$, and the family of maps $f_t$, $t \in [0,1]$ defined as the unique continuous solution of $h_t \circ f = f_t \circ h_t$, with $f_1=g$, gives an homotopy through Thurston maps. The desired isotopy follows. 
\gap

In the remainder of this work, we assume we are given a degree $d \geq 2$ Thurston map $f: S^2 \to S^2$ with postcritical set $P$. 
\gap 

\noindent{\bf Outline.} In \S 2, we collect generalities on Thurston maps and related combinatorial objects, 
e.g. the pullback relations on curves and on arcs. In \S 3, we discuss algebraic invariants and define 
the selfsimilarity complex, $\Sigma$. In \S 4, we derive dynamical consequences from the contraction 
of the virtual endomorphism.  \S 5 introduces the complex $\Gamma$ associated to an open cover $\UUU_0$, 
while \S 6 gives the proof that $\Sigma$ and $\Gamma$ can be chosen so that the natural map between them is a quasi-isometry. 
\S\S 7 and 8 conclude the proof of Theorem \ref{thm:characterization}. 
\gap 

\noindent{\bf Remark.} Nekrashevych \cite[Cor. 5.13]{nekrashevych:expanding} has shown that in great generality, a pair of maps $\iota, f: \MMM_1 \to \MMM_0$ satisfying analogous algebraic contraction conditions is combinatorially equivalent to a pair $\iota', f': \MMM_1' \to \MMM_0'$ where the $\MMM_i'$ are simplicial complexes, the maps $\iota', f'$ are piecewise affine, and the pair is expanding. While the result applies to our setting, the dimensions of $\MMM_i'$ can be larger than two, and it seems difficult to arrange for the spaces $\MMM_i'$ to coincide.

\section{Thurston maps}

\subsection{Curves and arcs} Here, we introduce the pullback relations on curves and arcs defined by a Thurston map, and formulate obstructions to expansion in terms of these relations. 

Denote by $\CCC$ the set of free homotopy classes of essential, {\em oriented}, simple, closed, nonperipheral (that is, not homotopic into arbitrarily small neighborhoods of elements of $P$) curves in $S^2-P$; we use the term {\em curve} for an element of $\CCC$. 

Let $o$ denote the union of the homotopy classes in $S^2-P$ of curves which are either inessential or peripheral; we call such curves {\em trivial}.  The {\em pullback relation} $\pullback$ on $\CCC\union \{o\}$ is defined by setting $o \pullback o$ and 
$ \gamma_1 \pullback \gamma_2$
if and only if  $\gamma_2$ is homotopic in $S^2-P$ to a component of $f^{-1}(\gamma_1)$.  Thus $\gamma \pullback o$ if and only if some preimage of $\gamma$ is inessential or peripheral in $S^2-P$.   
If $\gamma_1 \pullback \gamma_2$ we write $\gamma_1 \univpullback  \gamma_2$ if there is a representative of $\gamma_2$ mapping injectively to a representative of $\gamma_1$.   The relation defined by $\univpullback $ we refer to as the {\em univalent pullback relation}. 

We next define a similar pullback relation on certain arcs with endpoints in $P$.  Let $\AAA$ denote the set of isotopy classes (fixing $P$) of embedded
arcs $\alpha \subset S^2$ whose interiors are contained in $S^{2}-P$
and whose endpoints lie in $P$.   A {\em preimage} $\tilde{\alpha}$ of $\alpha$ is, by definition, the closure of a connected component of the inverse image of the interior of $\alpha$ under $f$.  If the endpoints of $\tilde{\alpha}$ lie in $P$, then $\tilde{\alpha}$ represents an element of $\AAA$.  Adjoining a symbol $o$ to the set $\AAA$ to stand for the collection of arcs with one or both endpoints outside $P$,  we obtain similarly a pullback relation on the set $\AAA\union \{o\}$.
\gap

We recall that a {\em multicurve} is a finite subset of $\CCC$ represented by disjoint curves.  A {\em Levy cycle} is a sequence $(\gamma_0, \gamma_1, \ldots, \gamma_{p-1})$ with $\gamma_i \univpullback  \gamma_{i+1 \ \bmod p}$ whose elements comprise a multicurve.  

The result below is not needed for the proofs of the theorems, since it could be derived from Theorems \ref{thm:characterization} and Theorem \ref{thm:prevents_levy} below.  We include the statement and proof since they are simple and help build intuition. It applies to maps which may have periodic branch points. 

\begin{prop}
\label{prop:obstructions_to_expansion}
Suppose $f$ is expanding. Given any curve $C_0$ or arc $\alpha_0$, there exists an integer $N$ such that if $C_0 \univpullback \ldots \univpullback C_n$ or $\alpha_0 \pullback \ldots \pullback \alpha_n$ is an orbit with $C_n, \alpha_n$ nontrivial, then $n \leq N$.  In particular,  the univalent pullback relation on curves has no (Levy) cycles or wandering curves, and the pullback relation on arcs has no cycles and no wandering arcs. 
\end{prop}

Abusing terminology somewhat, we summarize the conclusion by saying that the univalent pullback relation on curves and pullback relation on arcs have {\em no elements with arbitrarily long nontrivial iterates}. 
\gap

\pf  If $f$ were expanding with respect to a length structure, the conclusion would follow immediately, 
since lengths would decrease exponentially as one pulls back.   
Without this additional structure, we use instead coverings by small open balls. 
Let $\UUU_0, \UUU_i, i\in \N$ be as in the definition of expanding. 
Choose $m$ sufficiently large so that $C_0, \alpha_0$ are covered by some number, say $L$, of elements of $\UUU_m$ and the union of these elements are contained in a regular neighborhood of $C_0, \alpha_0$ which meets $P$ in the same manner as does $C_0, \alpha_0$ (that is, either not at all, or only at the endpoints of $\alpha_0$) and which is homotopic, relative to $P$, to $C_0, \alpha_0$.   It follows that for each $n$, $C_n, \alpha_n$ are covered by the same number $L$ of elements of $\UUU_{m+n}$.  Hence  if the diameters of the elements of the $\UUU_i$ tend to zero with $i$, the integer $n$ cannot be arbitrarily large. 
\qed

Since the pullback relations are natural with respect to combinatorial equivalence, Proposition \ref{prop:obstructions_to_expansion} gives a necessary condition for a Thurston map $f$ to be equivalent to an expanding map.  It is not, however, a sufficient condition; see \S \ref{secn:consequences} for examples.

\subsection{Orbifolds}  We review here the definition of the orbifolds associated to a Thurston map and their fundamental groups. 

For $x, y \in S^2$ let $\nu_0(y)=\lcm\{ \deg(f^n, x) : f^n(x)=y\}$
and $\nu_1(x)=\nu_0(f(x))/\deg(f,x)$.  For $i=0,1$ let $\Sigma_i=\{x \in S^2 : \nu_i(x)>1\}$, and let $\OOO_i$ be the orbifold
whose underlying topological space is $\{ x \in S^2: \nu_i(x)<\infty\}$ and whose weight function is $\nu_i$.  The sets $\Sigma_i$ are called the {\em singular sets} of $\OOO_i$.  Note that $\Sigma_0=P$. There are no singular points of infinite weight if and only if there are no periodic branch points of $f$.  In this case, the orbifolds $\OOO_i$ are compact.  
The singular sets satisfy  $f^{-1}(\Sigma_0) \supset \Sigma_1$.  Set $U_0 = S^2-\Sigma_0$ and $U_1=S^2-f^{-1}(\Sigma_0)$.  Then  $U_1 \subset U_0$ and $f: U_1 \to U_0$ is a covering
map.  

Let $b_0 \in U_0$ be a basepoint and $b_1\in f^{-1}(b_0)$ be one of
its preimages.  
For $i=0, 1$ let $N_i$ denote the normal subgroup of $\pi_1(U_i,
b_i)$ generated by the set of elements of the form $g^k$, where $g$
is represented by a simple closed peripheral loop $\gamma$
surrounding a puncture $x$ of $U_i$,  and the exponent $k$ is the
weight $\nu_i(x)<\infty$; if the weight is infinite, we do not add
such a loop as a generator.    The {\em orbifold fundamental groups}
$\pi_1(\OOO_i, b_i)$ are by definition the quotient groups
$\pi_1(U_i, b_i)/N_i$.   

In what follows, we put $\OOO=\OOO_0$ and $G:=\pi_1(\OOO, b_0)$. 

As the (orbifold) universal covering deck group, $G$ acts properly discontinuously on the universal covering of $\OOO$.

Note that each element of $\CCC$ corresponds to a conjugacy class in $G$ whose elements are of infinite order. 

\section{Algebraic constructions}

In this section, we briefly summarize constructions and results of Nekrashevych \cite{nekrashevych:book:selfsimilar}, specializing to the case of Thurston  maps. 

\subsection{The virtual endomorphism}
In the setup of \S 2.2, let $f_*: \pi_1(U_1, b_1) \to \pi_1(U_0, b_0)$ be the injective
homomorphism induced by the covering $f: U_1 \to U_0$.  Since $f$
sends peripheral loops to peripheral loops, it follows from the
definitions of the weight functions $\nu_i$ that $f_*: N_1 \to N_0$
is an isomorphism.  This observation and the ``Five Lemma" of
homological algebra imply that the homomorphism $f_*$ descends to a
well-defined and injective map $\cl{f}_*: \pi_1(\OOO_1, b_1) \to
\pi_1(\OOO_0, b_0)=G$.  We denote the image group
$\cl{f}_*(\pi_1(\OOO_1, b_1))$ by $H$; it has finite index in $G$.  

Let $\alpha: [0,1] \to U_0$ be a path joining $b_1$ to $b_0$ and
$\alpha_*: \pi_1(U_0, b_1) \to \pi_1(U_0, b_0)$ the induced
isomorphism.  Let $N_0'=\alpha_*^{-1}(N_0)$.   Since $N_0$ is normal,
the subgroup $N_0'$ is normal and is independent of the choice of
path $\alpha$.   Set $\pi_1(\OOO_0, b_1)=\pi_1(U_0, b_1)/N_0'$. 
Again, the map $\alpha_*$  descends to a well-defined isomorphism
$\cl{\alpha}_*: \pi_1(\OOO_0, b_1) \to \pi_1(\OOO_0, b_0)$.  

Since the inclusion $\iota: U_1 \hookrightarrow U_0$ sends peripheral
loops to loops which are either peripheral or trivial, and since
$\nu_0(x)$ divides $\nu_1(x)$ for all $x$, the induced map $\iota_*:
\pi_1(U_1, b_1) \to \pi_1(U_0, b_1)$ is surjective and sends $N_1$ to
$N_0'$.  It easily follows that the map $\iota_*$ also descends to a
surjective map $\cl{\iota}_*: \pi_1(\OOO_1, b_1) \to \pi_1(\OOO_0,
b_1)$.

\begin{defn}
\label{defn:virtual_endo}
The {\em virtual endomorphism induced by $f$} is the homomorphism
$\phi: H \to G$ defined by 
\[ \phi= \cl{\alpha}_* \circ \cl{\iota}_* \circ (\cl{f}_*)^{-1}.\]
\end{defn}
By construction, the virtual endomorphism $\phi$ associated to $f$ is
surjective.  

The virtual endomorphism depends on the choices of the basepoint $b_0$, the 
preimage $b_1$, and the homotopy class of the path $\alpha$.  Different
choices yield virtual endomorphisms which differ by pre- and/or
post-composition by inner automorphisms. Up to this ambiguity, 
the virtual endomorphism is an invariant of the homotopy class of $f$ relative to $P$. 
A combinatorial equivalence between Thurston maps conjugates, up to inner automorphisms, their virtual endomorphisms; the property of being contracting is preserved.

For $n \geq 2$, the {\em $n$th iterate} $\phi^n$ is the homomorphism whose domain  is defined
inductively by 
\[ \mbox{\rm dom}\phi = H; \;\;\; \mbox{\rm dom}  \phi^n=\{ g \in H :
\phi(g) \in \mbox{\rm dom}\phi^{n-1}\}\]
and whose rule is given by iterating $\phi$ a total of $n$ times.  

Suppose $S$ is a finite generating set for $G$.  We denote by $||g||$
the word length of $g$ in the generators $S$.
\begin{defn}
\label{defn:contracting1}
The virtual endomorphism $\phi: H \to G$ is called {\em contracting}
if the {\em contraction ratio}
\[ \rho= \limsup_{n\to\infty}\left( \limsup_{||g||\to\infty}
\frac{||\phi^n(g)||}{||g||}\right)^{1/n} < 1.\]
\end{defn}
The contraction ratio of the virtual endomorphism $\phi$ is
independent of the choices used in its construction.   

As an algebraic object, the virtual endomorphism is rather straightforward to describe. 
The next subsection gives a less familiar, but more natural, construction.  

\subsection{The biset $\frak{M}$}  

Suppose $\ell_0, \ell_1: [0,1] \to S^2 - P$ are two arcs joining the basepoint $b_0$ to a common point in $S^2 - P$.   
We say $\ell_0, \ell_1$ are {\em homotopic in $\OOO$} if the loop formed by traversing $\ell_0$ first and then the reverse of $\ell_1$ represents 
the trivial element of $G$. We denote by $\frak{M}$ the set of all homotopy classes $[\ell]$ of arcs $\ell$ in $\OOO$ joining
$b_0$ to one of the $d$ elements of $f^{-1}(b_0)$.  For $[\ell] \in \frak{M}$ we denote by $z_{[\ell]} \in f^{-1}(b_0)$ the corresponding common endpoint. 
The group $G$ acts on $\frak{M}$ via two commuting actions as follows. On the right, it acts by pre-concatenation:
\[ \ell \cdot g := \ell * g \]
where first $g$, then $\ell$ is traversed.  Note that the right action is free and has $d$ orbits.  The left action is obtained by taking preimages of the loops under $f$: 
\[ g \cdot \ell := f^{-1}(g)[z_{[\ell]}]* \ell \]
where first $\ell$, then the lift $f^{-1}(g)[z_{[\ell]}]$ of $g$ based at $z_{[\ell]}$, is traversed.  
The set $\frak{M}$ equipped with these two commuting actions of $G$ defines the $G$-{\em biset} associated to $f$. 

In the next subsection, we show that choosing certain representatives of elements of $\frak{M}$ and representatives of generators of $G$ gives rise to an abstract $1$-complex $\Sigma$ and a projection map $\pi_\Sigma: \Sigma \to S^2$. Along the way, we reformulate the notion of contraction in terms of an action of $G$ on the set of words in an alphabet $X$ of size $d$. 

\subsection{The selfsimilarity complex $\Sigma$}

Choose a bijection $\Lambda: X \to f^{-1}(b_0)$, where $X$ is a finite set of cardinality $d$. 
For each $x \in X$, let $\lambda_x$ be an arc in $S^2 - P_f$ joining $b_0$ to $\Lambda(x)$.
The collection $\{[\lambda_x]: x \in X\}$ of elements of $\frak{M}$ is called a {\em basis} of the biset $\frak{M}=\frak{M}_f$. 

For $n \in \N$, consider now the $n$th iterate $f^n$ of $f$. The orbifolds $\OOO=\OOO_f$ and $\OOO_{f^n}$ coincide, so $G=G_f = G_{f^n}$.  By lifting the arcs $\lambda_x$ under $f^n$, we obtain an identification $\Lambda: X^n \to f^{-n}(b_0) \times \{n\}$ of the set $X^n$ of words of length $n$ in the alphabet $X$ with the fiber $f^{-n}(b_0) \times \{n\}$ and a corresponding basis for the $G$-biset $\frak{M}_{f^n}$ associated to the $n$th iterate of $f$.   We denote by $X^* = \union_n X^n$, with $X_0 = \{\emptyset\}$ consisting of the empty word.   We write $|w|=n$ if $w \in X^n$. 

Suppose $S$ is a generating set for $G$.  For each $s \in S$, let $\gamma_s$ be a loop in $S^2-P$ based at $b_0$ representing $s$.  By lifting loops under iterates of $f$, we obtain an action of $G$ on $X^*$ preserving the lengths of words and acting transitively on each $X^n$. 

The basepoint $b_0$, arcs $\lambda_x, x \in X$, and loops $\gamma_s, s \in S$, we refer to as {\em defining data}.

Since the right action is free, for each word $u \in X^*$, and each $g \in G$, there are a unique word $v \in X^n$ and a unique $h \in G$ with $g.u = v.h$; the element $h$ is called the {\em restriction}\footnote{sometimes the term ``section'' is used.} of $g$ to $u$, and is denoted $g|_u$.  The interpretation in terms of defining data is as follows. Let $\lambda_u, \lambda_v$ denote the concatenations of lifts of arcs corresponding to $u$ and $v$, so that the endpoints of $\lambda_u, \lambda_v$ correspond to $u$ and $v$.  Then $v$ is the endpoint of $f^{-n}(g)[u]$ and $h$ is represented by the loop which first traverses $\lambda_u$, then $f^{-n}(g)[u]$, then the reverse of $\lambda_v$.  

By \cite[Prop. 2.11.11]{nekrashevych:book:selfsimilar}, contraction of the virtual endomorphism is equivalent to the following: {\em there is a finite set $\NNN$ such that for each $g \in G$, there is an integer $n$ such that for all $v \in X^*$ satisfying $|v| \geq n$, the restriction $g|_v \in \NNN$}.  In this case, one says the biset $\frak{M}$ is {\em hyperbolic}.  This is an asymptotic condition: $\frak{M}_f$ is hyperbolic if and only if $\frak{M}_{f^m}$ is hyperbolic for some $m \geq 1$, since the contraction coefficients for $\phi$ and $\phi^m$ are necessarily simultaneously less than $1$. 

Suppose now that defining data has been chosen. The {\em selfsimilarity complex} $\Sigma$ is the infinite $1$-complex whose vertex set is $X^*$ and whose edges, which come in two types, {\em horizontal} and {\em vertical}, are defined as follows. Given $w \in X^*$, a horizontal edge, labelled $s$, joins $w$ to $s.w$; a vertical edge, labelled $x$, joins $w$ to $xw$.  The selfsimilarity complex, with edges of length $1$, is a proper geodesic metric space, and the local valence at each nonempty vertex is the same if $S$ is symmetric. 

The right shift $\sigma: X^* \to X^*$ determines a cellular self-map $f_\Sigma: \Sigma' \to \Sigma$ where $\Sigma'$ is the subcomplex determined by nonempty words.  When $\frak{M}$ is hyperbolic, $\Sigma$ is hyperbolic in the sense of Gromov, and the boundary $\JJJ$ is compact.  Since the virtual endomorphism 
$\phi$ is surjective, the action is recurrent, so $\JJJ$ is also connected and locally connected \cite[Theorem 3.5.1]{nekrashevych:book:selfsimilar}.  The map $f_\Sigma$ extends to a map of the boundary $f_\Sigma: \JJJ \to \JJJ$, yielding a dynamical system. 

Although $\Sigma$ is an abstract $1$-complex, path-lifting of the defining data determines a map 
$\pi_\Sigma: \Sigma \to S^2$; see \cite{kmp:gromov} and \cite[Theorem 5.5.3]{nekrashevych:book:selfsimilar}.  
If $f$ is not expanding with respect to a length structure on $S^2$, however, $\pi_\Sigma$ need not 
extend to a map $\JJJ \to S^2$.

\section{Dynamical consequences of the hyperbolicity of $\frak{M}$}
\label{secn:consequences}

Suppose $f$ has no periodic branch points. Proposition \ref{prop:obstructions_to_expansion}  gives necessary conditions for $f$ to be expanding. The main result of this section is to derive the same conclusion from the algebraic assumption that $\phi$ is contracting; equivalently, that $\frak{M}$ is hyperbolic.

\begin{thm}
\label{thm:prevents_levy}
Suppose $f$ has no periodic branch points. If $\frak{M}$ is hyperbolic, the univalent pullback relation  on curves and the pullback relation on arcs have no elements with arbitrarily long nontrivial iterates.
\end{thm}

\pf The central ingredient in the proof is the following Lemma, which will be proved later.  
\begin{lemma}
\label{lemma:topology_algebra}
Suppose $C_0\univpullback C_1$.  Given any $g_0 \in G$ whose conjugacy class corresponds to $C_0$, there exists $g_1 \in G$ whose conjugacy class corresponds to $C_1$ and $x \in X$ with $g_0 \cdot x = x \cdot g_1 \in \frak{M}$.  
\end{lemma}
To prove the Theorem, suppose $C_0  \univpullback C_1 \univpullback \ldots \univpullback C_n$.  Induction and Lemma \ref{lemma:topology_algebra} implies that there exist $g_i \in G, x_i \in X$ such that 
\[ g_0 \cdot x_0x_1\ldots x_i = x_0x_1\ldots x_i \cdot g_{i+1}, \;\;\; i=0, \ldots, n-1.\]
In particular, upon setting $w=x_0x_1 \ldots x_n \in X^{n+1}$, we have 
\[ g_0 \cdot w = w \cdot g_n;\]
recall that by definition, $g_n = g_0|_w$. 
Since $\frak{M}$ is hyperbolic, there is a finite set $\NNN \subset G$ such that $g_n =g_0|_w \in \NNN$ whenever $n=|w|$ is sufficiently large.  
This immediately implies that the number of distinct elements $g_n$ of $\CCC$ arising in this way is bounded and hence that the collection of curves $\{C_n\}$ arising in this way is also bounded, independent of $n$. 

To rule out cycles, we use topology to simplify the algebra.  Consider first the case of a Levy cycle of period one, i.e. a curve $C$ for which $C\univpullback C$.  By varying $f$ within its homotopy class, and using the fact that $C$ is oriented, we may assume $f=\id$ on a neighborhood of $C$.  Take the basepoint $b_0$ for $\frak{M}$ to lie on $C$; it becomes a fixed point of $f$.  Let $g$ be represented by $C$, regarded as a loop based at $b_0$.  By choosing a basis for $\frak{M}$ with one arc equal to the constant path at $b_0$ and corresponding to an element $x$ in the alphabet $X$, we find $g \cdot x = x \cdot g$.  Induction shows then that for all $n \in \N$ and $m \in \Z$ we have that if $w=xx\ldots x \in X^n$ then 
\[ g^m \cdot w = w \cdot g^m.\]
But this is impossible if $\frak{M}$ is hyperbolic:  $g$ has infinite order,  and the identity above shows that for all $n$, $g^m|_w = g^m$ which cannot lie in $\NNN$ if $m$ is sufficiently large.   

The case of a Levy cycle of period larger than one can be ruled out by passing to an iterate.   

Analyzing the pullback relation on arcs is reduced to that of the pullback relation on curves.  An arc which does not eventually become trivial must have both endpoints in cycles of $P$. For such arcs $\alpha_0, \alpha_1$, since $f$ has no periodic branch points, whenever $\alpha_0 \pullback \alpha_1$, the boundaries of  regular neighborhoods of $\alpha_0$ and $\alpha_1$ yield curves $C_0$ and $C_1$ for which $C_0 \univpullback C_1$. 

If $\#P> 3$ then $C_0, C_1$ must be nontrivial, and the eventual triviality of curve $\univpullback$-orbits, already established, implies eventual triviality of arc $\pullback$-orbits. So assume $\#P = 3$. By Thurston's characterization \cite{DH1}, $f$ is equivalent to a rational function without periodic critical points. Such a map is uniformly expanding with respect to a complete length structure, so such arc orbits are always eventually trivial. 
\qed

We now turn to the proof of the Lemma, which is just an exercise in definitions.
\gap

\pf Suppose  defining data has been chosen. Choose a representative $\gamma_0$ for $g$, and let 
$\wtgamma_0$ be a preimage of $\gamma_0$ corresponding to $C_1$; it exists, by homotopy lifting and the fact that 
$C_1$ maps to $C_0$ by degree one.   Let $\tilde{b}$ be the unique preimage of $b$ lying on $\wtgamma_0$, 
and suppose $x$ corresponds to $\tilde{b}$, i.e. $\Lambda(x)=\tilde{b}$.  The loop given by 
\[ \gamma_1 := \cl{\lambda}_x*\wtgamma_0 * \lambda_x\]
(where $\lambda_x$ is traversed first) represents an element $g_1 \in G$.  Moreover, by the definition of the biset, 
\[ x \cdot g_1 = [\lambda_x * \cl{\lambda}_x * \wtgamma_0 * \lambda] = g_0 \cdot x.\]
\qedspecial{Lemma \ref{lemma:topology_algebra}}

Examples showing the conditions in Propositions \ref{prop:obstructions_to_expansion} and Theorem \ref{thm:prevents_levy} are necessary but not sufficient can be found among Thurston maps induced by  torus endomorphisms that are hyperbolic but not expanding. 

Let $T^2 = \R^2/\Z^2$, let $A$ be an integral matrix with determinant larger than one, and let $F: T^2 \to T^2$ be the noninvertible endomorphism on the torus induced by $A$. There is a corresponding pullback relation on curves, defined analogously as in \S 2. In this setting, all preimages of nontrivial curves are nontrivial, are homotopic to each other, and map by the same degree.   

\begin{prop}
\label{prop:eval_condition}
If a curve $C$ has arbitrarily long nontrivial iterates under the univalent pullback relation, then either $C \univpullback C$ or $C \univpullback -C$, where $-C$ denotes the curve $C$ with the opposite orientation.
This occurs if and only if  the matrix $A$ has an eigenvalue $\lambda$ equal to $+1$ or $-1$. 
\end{prop}

\pf The sufficiency in the second statement is obvious. The proof will show both the first statement, and the necessity; it is a straightforward reinterpretation of the arguments of \cite[Prop. 2.9.2]{nekrashevych:book:selfsimilar}. 
Use the standard basis to identify $\Q^2$ with $H_1(T^2, \Q)$ and $\pi_1(T^2)$ with $H_1(T^2, \Z)=\Z^2$, so that $A$ is the matrix for the induced map on rational and integral homology. Let $B=A^{-1}$. Observe that if $C=C_0 \univpullback C_1$ then in fact all preimages of $C_0$ map by degree one and the monodromy corresponding to $C_0$ is trivial. Thus if $C_0 \univpullback C_1 \univpullback C_2 \univpullback \ldots $ is an infinite sequence of related curves (either wandering, or cycling) then as an element of $\pi_1(T^2)$, the curve $C_0$ represents a nontrivial element of the kernel $U < \Z^2< \Q^2$ of the corresponding iterated monodromy group action.  The group $U$ is $B$-invariant, so $U \otimes \Q$ is a $B$-invariant $\Q$-linear subspace of $\Q^2$; let $C=B|_{U \otimes \Q}$. Then the characteristic polynomial $q$ of $C$ is monic, has integer coefficients, and is a factor of the characteristic polynomial of $B$. If $1/\lambda$ is a root of $q$ then $\lambda$ is an eigenvalue of $A$. Both $\lambda$ and $1/\lambda$ are algebraic integers, hence the norm of $\lambda$ is equal to $1$. The condition $\det(A)>1$ implies the algebraic degree of $\lambda$ must be equal to $1$, so $\lambda = \pm 1$. 
\qed

The endomorphism $F: T^2 \to T^2$ is {\em hyperbolic} if $A$ has no eigenvalues of unit modulus.  Under the assumption $\det(A)>1$, the following is easily shown via straightforward computations. 

\begin{prop}
\label{prop:hyperbolic_endo}
The endomorphism $F$ is hyperbolic but not expanding if and only if $\trace(A)-\det(A)> 1$. 
\end{prop}

The endomorphism $F: T^2 \to T^2$ commutes with $\pm \id_{T^2}$ and descends to a Thurston  map $f: S^2 \to S^2$; the corresponding orbifold $\OOO$ has signature $(2,2,2,2)$. 
The singular points are the images of the points of order at most $2$. 

\begin{thm}
\label{thm:not_expanding}
Suppose $\det(A)>1$ and $\trace(A)-\det(A)> 1$, and let $f: S^2 \to S^2$ be the induced Thurston map. Then 
\be
\item $f$ is not equivalent to an expanding map, but 
\item the univalent pullback relation on curves and the pullback relation on arcs have no elements with arbitrarily long nontrivial iterates. 
\eb
\end{thm}

\pf Proposition, \ref{prop:hyperbolic_endo}, Proposition \ref{prop:eval_condition}, and Theorem \ref{thm:torusversion} imply both conclusion (1) and the first part of (2).  Suppose $\alpha \subset S^2$ is an arc. The two preimages of $\alpha$ under the branched covering $T^2 \to S^2$, when concatenated, yield a nontrivial curve, $C$.  If $\alpha_0 \pullback \alpha_1$ then for the corresponding curves in $T^2$ we have $C_0 \univpullback C_1$. Thus nontrivial arc orbits lift to nontrivial univalent curve orbits.
\qed

\section{The covering complex $\mathbf{\Gamma}$} \label{scn:covering}

In this section, we begin by summarizing constructions and results of \cite[Chapter 3]{kmp:ph:cxci}, 
specializing to the case of Thurston  maps. We associate to  $f$ and a suitable open covering $\UUU_0$ another infinite 
one-complex, $\Gamma$, this time a graph. When it is hyperbolic, it has a natural boundary, and there is 
a natural map $\pi_\Gamma: S^2 \to \bdry \Gamma$. 

For convenience, equip $S^2$ with the usual Euclidean length metric. 
Let $\UUU_0$ be a finite open cover of $S^2$ by small balls $B$ for which $\#B \intersect P \leq 1$ 
with equality only when the corresponding point of $P$ is in the center of $B$.   
For $n \in \N$ let $\UUU_n$ be the open cover whose elements are connected components of sets of the form 
$f^{-n}(U), U \in \UUU_0$. It is convenient to set $\UUU_{-1}:=\{S^2\}$.  
If $U \in \UUU_n$ we write $|U|=n$. A key observation here is that if $f$ has no periodic branch points, 
then there is an upper bound on the degrees of restrictions $f^n: \wtU \to U$, $\wtU \in \UUU_n, n \in \N$.

The vertex set of $\Gamma$ is defined as $\union_{n\geq -1}\UUU_n$, i.e. one vertex for each component.   
We denote the vertex corresponding to the unique element of $\UUU_{-1}$ by $o$.  
Horizontal edges join $U_1, U_2$ if $|U_1|=|U_2|$ and $U_1 \intersect U_2 \neq \emptyset$. Vertical edges join $U, V$ if $||U|-|V|| =1$ and $U \intersect V \neq \emptyset$.   $\Gamma$ is a graph with root $o$.  We equip $\Gamma$ with the graph path metric, so that each edge has length $1$. 
Then $\Gamma$ is a proper, geodesic metric space; let $d$ denote the corresponding distance function. 
If $f$ has no periodic branch points, then the valence of a vertex of $\Gamma$ is uniformly bounded, 
independent of the vertex. 

Let $\varepsilon>0$ be small, and let $\varrho_\varepsilon: \Gamma \to (0,\infty)$ be given by 
$\varrho_\varepsilon(x)=\exp(-\varepsilon d(x, o))$.   Upon redefining the lengths of curves via 
\[ \ell_\varepsilon(\gamma):=\int_\gamma  \varrho_\varepsilon \, ds\] 
we obtain an incomplete length metric space. Its completion compactifies $\Gamma$
by adding a boundary at infinity $\partial_{\ep}\Gamma$. 
Points in its boundary $\bdry_\ep \Gamma$ are represented by vertical geodesic rays from the root.  Given $a \in S^2$, let $U_n(a), n=-1, 0, 1, 2, \ldots$ be a sequence of open sets for which $a \in U_n \in \UUU_n$. Then $\{U_n\}$ defines a geodesic ray and, therefore, a point of $\bdry_\ep \Gamma$. 
We obtain a well-defined surjective continuous map $\pi_\Gamma: S^2 \to \bdry_\ep \Gamma$. In particular, $\bdry_\ep \Gamma$ is connected 
and locally connected. When $\Gamma$ is hyperbolic and $\ep>0$ is chosen small enough, 
then $\bdry_\ep\Gamma$ coincides with its visual boundary $\bdry\Gamma$. We will always choose $\ep$ as above when $\Gamma$ is hyperbolic.

\begin{lemma}
\label{lemma:monotone}
If $\Gamma$ is hyperbolic, the projection $\pi_\Gamma: S^2 \to \bdry\Gamma$ is monotone, i.e. its fibers are connected.
\end{lemma}

\pf For $z\in S^2$ and $n\ge -1$, pick $U_n(z)\in\UUU_n$ such that $z\in U_n(z)$. It follows that the sequence $(U_n(z))_n$ defines a ray in $\Gamma$.
Suppose $\pi_\Gamma(x)=\pi_\Gamma(y)=\xi \in \bdry \Gamma$.  From the definition of $\pi$ and the hyperbolicity of $\Gamma$, 
there is a positive integer $D$ such that for every $n \in \N$, there exist a ``chain'' 
$U_n^0, U_n^1, \ldots, U_n^D$ of elements of $\UUU_n$ 
such that $x\in U_n^0$, $y \in U_n^D$, and $U_n^{i-1}\intersect U_n^i \neq \emptyset$ for all $i=1, \ldots, D$.  
For $n \in \N$ let  $Y_n = U_n^0 \union \ldots\union U_n^D$.  Then $Y_1, Y_2, Y_3, \ldots$ is a sequence of open connected sets containing both $x$ and $y$. 
Set $$K=\cap_{k\ge 0} \overline{\cup_{n\ge k} Y_n}\,.$$  It follows that  $K\subset S^2$ is a continuum containing both $x$ and $y$.
We will prove that $\pi_\Gamma(K)=\xi$ and, therefore, that $\pi_\Gamma^{-1}(\xi)$ is connected. 
Fix $z \in K$.  For any $k\ge 0$, $U_k(z)\cap K\ne\emptyset$, so we may find $n\ge k$ such that $Y_n\cap U_k(z)\ne\emptyset$.
This means that 
\begin{eqnarray*}
d_{\ep}(\pi_\Gamma(z),\xi)& \le & d_{\ep}(\pi_{\Gamma}(z),U_k(z))+ d_{\ep}(U_k(z),\{U_n^j, 0\le j\le D\})\\
& & \quad +\diam_{\ep} \{U_n^j, 0\le j\le D\} + d_{\ep}(\xi,\{U_n^j, 0\le j\le D\})\\
& \lesssim & e^{-\ep k} + e^{-\ep k} + e^{-\ep n} + e^{-\ep n}\lesssim e^{-\ep k}\,.\end{eqnarray*}
Since $k$ is arbitrary, we obtain that $\pi_{\Gamma}(z)=\xi$.
\qed

From \cite[Theorem 5]{daverman:decompositions}, we have that the decomposition of $S^2$ by the fibers of $\pi_\Gamma$ is upper-semicontinuous, and the quotient space is homeomorphic, via the map induced by $\pi_\Gamma$, to $\bdry\Gamma$.

The map $f$ induces a cellular map $f_\Gamma: \Gamma' \to \Gamma$, where $\Gamma'$ is the subgraph on vertices at level $\geq 1$.  Since it sends vertical rays to vertical rays,  it induces a continuous map $f_\Gamma: \bdry \Gamma \to \bdry \Gamma$.
Since it sends sufficiently small round balls to round balls \cite[Prop. 3.2.2]{kmp:ph:cxci}, it is an open map; it is easily seen to be at most $d$-to-$1$, so it is closed as well.  In particular, the local degree function 
\[ \deg(f_\Gamma,\xi):=\inf_W\{\#{f_\Gamma}^{-1}(\zeta) \intersect W : \zeta \in f_\Gamma(W), \;\; \xi \in W \; \mbox{open} \}\]
is well-defined.

\section{Coincidence of quasi-isometry types}
\label{secn:coincidence} 

Recall that if $X, Y$ are metric spaces, a function $\Phi: X \to Y$ is a {\em quasi-isometry} if there exist constants $\lambda \geq 1, c \geq 0, b \geq 0$ such that 
\[ \frac{1}{\lambda}|x-x'| - c \leq |\Phi(x)-\Phi(x')| \leq \lambda|x-x'| + c \]
and such that $|y - \Phi(X)|\leq b$ for all $y \in Y$; here $|x-x'|=d_X(x,x')$, etc..  For geodesic metric spaces, hyperbolicity is preserved under quasi-isometries.

Recall that a choice of defining data defines a projection $\pi_\Sigma: \Sigma \to S^2$ and an identification map $\Lambda:X^{*}\to \cup_{n\ge 0} f^{-n}(\{b\})$. The main result of this section is 

\begin{thm}
\label{thm:coincidence}
Let $f: S^2 \to S^2$ be a Thurston map without periodic critical points.
Then there exist defining data for the associated selfsimilarity complex $\Sigma$ and the covering graph $\Gamma$ 
and a  level-preserving quasi-isometry $\Phi: \Sigma \to \Gamma$ such that,
for every vertex $v\in\Sigma^0$, we have $\Lambda(v)\in\Phi(v)$ and  $\Phi \circ f_\Sigma = f_\Gamma \circ \Phi$.
\end{thm}

\begin{remark} If we admit periodic branch points, then the same proof will show that $\Sigma$
is quasi-isometric to a covering complex $\Gamma$ built from a cover of the complement in the sphere of a neighborhood of the set of periodic branch points.
\end{remark}

\pf We first choose the defining data  and list some elementary properties.  
For convenience, we endow $S^2$ with its ordinary Euclidean length metric.  

\gap

\noindent{\bf Choosing the defining data.}  
Choose the basepoint $b \in S^2-P$ to be a regular point.  We may choose representative loops 
$\gamma_s, s \in S$, so that each  is simple, peripheral, and has the property that for some common $\epsilon>0$, 
the closed $\epsilon$-neighborhood $N_\epsilon(\gamma_s)$ is contained in a closed annulus $A$ contained in $S^2-P$ surrounding on one side 
at most one point of $P$.  The peripheral condition and the assumption that there are no periodic critical points implies that 
there exists $N_0$ such that for all $n$, 
\[ \deg(f^n: \widetilde{A} \to A) \leq N_0\] 
for all preimages $\widetilde{A}$ of $A$ under $f^{-n}$.  
In this situation, if the above degree is say $\delta$, then in the graph $\Lambda$, for each vertex $v$ at level $n$ with $\Lambda(v) \in \widetilde{A}$, there is a closed horizontal edge-
path of length $\delta$ whose edges are labelled by $s$ corresponding and which corresponds to the monodromy at level $n$ induced by $\gamma_s$.

We may also choose representative arcs $\lambda_x \subset S^2-P, x \in X$, 
so that each is simple and has the property  that for some common $\epsilon>0$, 
the closed $\epsilon$-neighborhood $N_\epsilon(\lambda_x)$ is a closed disk $D$ contained in $S^2-P$.  Thus
\[ \deg(f^n:\widetilde{D} \to D) = 1\] 
for all preimages $\widetilde{D}$ of $D$ under $f^{-n}$.

For $w\in X^*$, let $\lambda_{xw}$ be the lift of $\lambda_{x}$ under $f^{|w|}$ which starts at $\pi_{\Sigma}(w)$.   

\gap

Since $f$ is a branched covering, given $\delta>0$, there is a constant $r(\delta)>0$ with $r(\delta) \to 0$ as $\delta \to 0$ 
such that for any ball $B(y,\delta) \subset S^2$, any preimage of $B$ under $f$ has diameter $<r(\delta)$.  
It follows that we may choose the defining data $\UUU_0$ for the covering graph $\Gamma$  to consist of a finite collection 
of small open balls satisfying the following properties.

Given a set $E \subset S^2$, we denote by 
\[ \UUU_n(E) = \{ U \in \UUU_n : E \intersect U \neq \emptyset\}.\]

\be
\item $\forall U \in \UUU_0$, $P \intersect \bdry U = \emptyset$, $\# P \intersect U \leq 1$, and if equality holds, then the postcritical point in $U$ is at the center of $U$; 
\item $\forall U \in \UUU_0$, $\UUU_0(U) \intersect P \leq 1$;
\item $\forall x \in X$, $\cl{\UUU_0(\lambda_x)} \subset N_\epsilon(\lambda_x)$;
\item $\forall s \in S$, $\cl{\UUU_0(\gamma_s)} \subset N_\epsilon(\gamma_s)$;
\item $\forall U_0 \in \UUU_0$ and $U_1 \in \UUU_1$ with $U_0 \intersect U_1 \neq \emptyset$, 
the union $U_0 \union U_1$ is contained in a ball $B$ such that the ball $2B$ of twice the radius meets at most one point of $P$; 
\item $\forall p \in S^2$ there exists an arc $\alpha$ joining the basepoint $b$ to $p$ such that $\UUU_0(\alpha)$ is contained in a disk which meets $P$ in at most one point.
\item There is exactly one element of $\UUU_0$ containing the basepoint $b$, and this basepoint is the center of the ball.  
\eb

{\noindent\bf Basic properties.} \\
\gap
\noindent{\bf (BP1).} Condition 1 implies that every element $U\in\UUU_n$  has a unique preferred point $c_U$
(we call it the center of $U$) which is the unique preimage of the center $c$ of the ball $f^{\circ n}(U)$ under $f^n$. 

\gap

{\noindent\bf (BP2).} Conditions 1 and 2 imply the following. Since there are no periodic branch points, 
there is a uniform upper bound $q$ on the local degree of all iterates of $f$.   
Suppose $f^n:\wtU_i \to U_i \in \UUU_0, i=0,1$.  If  $\wtU_0, \wtU_1 \in \UUU_n$ intersect, 
then at most one of $\wtU_0, \wtU_1$ meets a branch point of $f^n$; say it is $\wtU_0$.  
Then there are $q' \leq q-1$ elements $\wtU_2, \ldots, \wtU_{q'}$ of $\UUU_n$, each a preimage of $U_1$, 
such that $f^{\circ n}: \wtU_0 \union \wtU_1 \union \ldots \union \wtU_{q'} \to U_0 \union U_1$ is proper and of degree $q'$.

\gap

{\noindent\bf (BP3).} From condition 6, we define $A_0$ to be the finite set of curves $\alpha_U:[0,1]\to S^2$ joining the center 
$c_U =\alpha_U(0)$ of $U$
to $b=\alpha_U(1)$. Let also $A_n$ denote all the lifts of these curves under $f^n$: these are curves joining every center
of $U\in\UUU_n$ to a point of $f^{-n}(b)$. Since every $\alpha_U$ is an arc contained in a disk which intersects $P$
in at most one single point, there is some $L$ independent of $n$ and $U\in\UUU_n$ 
such that $\alpha_U$ is covered by at most $L$ elements of $\UUU_n$.

\gap

{\noindent\bf (BP4).} Condition 7  is merely to make the exposition in what follows more convenient, 
for with this choice we will define  $\Phi: \Sigma \to \Gamma$  in a simple way.  

\gap

\noindent{\bf Strategy of the proof. } Define $\Phi:V(\Sigma)\to V(\Gamma)$ by setting $\Phi(w)$ to be the unique set 
$U\in\UUU_{|w|}$ which
contains $\Lambda(w)$. We extend $\Phi$ over the edges to all of $\Sigma$ as follows. 
Suppose $\gamma_s$ represents a horizontal edge at level $0$ (a loop) in $\Sigma$. In the sphere, as this edge is traversed, 
it passes through a sequence of elements of $\UUU_0$; we pick such a sequence, and associate to this a horizontal loop $\Phi(\gamma_s)$ in $\Gamma$ at level $0$, 
parameterized at constant speed. We define $\Phi$ to send the edge of $\Sigma$ corresponding to $\gamma_s$ to this loop $\Phi(\gamma_s)$. 
Similarly, each horizontal edge $e=e_s$ at level $n$ in $\Sigma^1$ corresponds to
a lift under $f^n$ of a unique edge $\gamma_s$, $s\in S$, at level zero; hence by lifting the chosen covering for $\gamma_s$ we obtain again  
a curve $\Phi(e)$ in $\Gamma$ joining the images under $\Phi$ of
their extremities by definition.   

Similarly, for $x\in X$, we associate to the arc $\lambda_x$  an edge-path $\Phi(\lambda_x)$ in $\Gamma$ corresponding to a covering of the chosen arc $\alpha_{\Lambda(x)}$ by elements of $\UUU_0$ 
concatenated with the vertical edge of $\Gamma^1$ between the elements $U\in\UUU_0$ and $V\in\UUU_1 $ containing $\Lambda(x)$.
We again extend $\Phi$ to the other vertical edges 
equivariantly by path-lifting.

We also define a coarse inverse $\Psi : V(\Gamma)\to V(\Sigma)$ by setting $\Psi(U)=w$, where $w\in V(\Sigma)$ is the vertex such that
$\Lambda(w)$ is the endpoint of $\alpha_U$.  Note that $\Psi\circ\Phi$ is the identity on $V(\Sigma)$ since $\alpha_U$ 
is trivial for preimages
 of $U_b\in\UUU_0$ containing $b$, and $\Phi\circ\Psi$ maps a vertex $U\in\Gamma^0$ to the endpoint of $\alpha_U$.

\gap

To prove that $\Phi$ is a quasi-isometry, we will show that both $\Phi$ and $\Psi$ are Lipschitz on the set of vertices 
and that the image of $\Phi$ is cobounded. The conclusion will follow since  $\Psi\circ\Phi$ is the identity on $\Sigma^0$.
To prove that each map is Lipschitz, we will show that the distance between the images
of the two extremities of an edge is uniformly bounded.

\gap

\noindent{\bf $\mathbf{\Phi} $ is Lipschitz.}   
Let us consider the horizontal case first.  Suppose $v, w \in V(\Sigma)$, $|v|=|w|=n$, and $v, w$ are joined by a horizontal 
edge, so that $w=s.v$ for some $s \in S$.  Equivalently, there is a lift $\wtgamma_s$ of $\gamma_s$ under $f^n$ joining 
$\pi_\Sigma(v)$ to $\pi_\Sigma(w)$. Note that $f^n:\wtgamma_s\to \gamma_s$ is $1-1$ away from the extremities. 
The curve $\gamma_s$ is covered by $\#\UUU_0(\gamma_s)$ elements of $\UUU_0$.  The lift $\wtgamma_s$ is covered by at most $\#\UUU_0(\gamma_s)+1$ elements of $\UUU_n$.   
By the definition of the metric in $\Gamma$, this implies that $|\Phi(v)-\Phi(w)|_{hor} \leq  (\#\UUU_0(\gamma_s) +1)$.  Note that we may assume that $ \#\UUU_0(\gamma_s)\le
\#\UUU_0-1$.
Since $|\Phi(v)-\Phi(w)|_\Gamma \leq |\Phi(v)-\Phi(w)|_{hor}$, we conclude that 
\[ |\Phi(v)-\Phi(w)|_\Gamma \leq C_{hor}:= \#\UUU_0. \]

Now let us consider the vertical case.  Suppose $|v|=n+1$, $|w|=n$, and $v$ and $w$ are joined by a  vertical edge in $\Sigma$.  By definition, $w=xv$ for some $x \in X$, and the points $\Lambda(w)$ and $\Lambda(v)$ in $S^2$ are joined by a lift of the arc $\lambda_x$ under $f^{-n}$ starting at $\Lambda(v)$. 
The arc $\lambda_{x}$ is covered by $\#\UUU_0(\lambda_x)$ elements of $\UUU_0$, 
each contained in  the neighborhood $N_\epsilon(\lambda_x)$, which is a disk in $S^2-P$.  
It follows that $\widetilde{\lambda}_x$ is covered by the same number $\#\UUU_0(\lambda_x)$ of elements of $\UUU_n$.  
By the definition of the metric in $\Gamma$, this implies that $|\Phi(v)-\Phi(w)|_\Gamma \leq \#\UUU_0(\lambda_x) +1$.  Hence 
\[ |\Phi(v)-\Phi(w)|_\Gamma \leq C_{ver}:=\max\{\#\UUU_0(\lambda_x) : x \in X\}+1.\]

We conclude that $\Phi$ is $C=\max\{C_{hor}, C_{ver}\}$-Lipschitz on $V(\Sigma)$.  
\gap

\noindent{\bf The image of $\Phi$ is cobounded in $\Sigma$.}  The proof of this follows along exactly the same lines.   
For any $U\in \UUU_n$,  the curve $\alpha_U\in A_n $ joins $c_U$ to $\Phi\circ\Psi(U)$, whose
distance is bounded by $L$ according to the basic property (BP3).

\gap

\noindent{\bf Realization of the covering complex $\Gamma$.} Proving that $\Psi$ is Lipschitz 
will require some more data.  
The covering complex $\Gamma$ comes equipped with a self-map $f_{\Gamma}:\Gamma\setminus B(X,1)\to\Gamma$. 
We will find a continuous map $\Gamma \to S^2$ which conjugates $f_{\Gamma}$ and $f$ where defined.
First, if $e=(U_1,U_2)$ denotes an edge with $U_1,U_2\in \UUU_0$, we consider $\beta_e$ to be the geodesic
(in the round metric) joining their centers; recall that $U_j$ are spherical balls at level $0$, so that 
$\beta_e\subset U_1\cup U_2$. Denote by $B_0$ this finite set of curves, and let $B_n$ denote the lifts of the elements
of $B_0$ under $f^n$. 


If $e=(U_0,U_1)$ is a vertical edge of $\Gamma$ with $U_0\in\UUU_0$ and $U_1\in\UUU_1$, we may find an arc $\eta_e$
joining their centers inside $U_0\cup U_1$. Let $H_0$ denote the sets of arcs defined for every such edge, and for each $n \geq 0$, define
$H_n$ to be the set of their lifts under $f^n$. 

\gap

\noindent{\bf $\mathbf{\Psi:=\Phi^{-1}}$ is Lipschitz.}   Let us first consider an edge $e$ joining two sets 
$\wtU_1$, $\wtU_2$ at the same level $n\ge 0$.
Let $\widetilde{\gamma}=\alpha_{\wtU_1} * \beta_e *\overline{\alpha_{\wtU_2}}$. 
The image $f^n(\widetilde{\gamma})$ is a closed curve based at $b$ which is homotopically nontrivial 
in $S^2\setminus P$ since $\widetilde{\gamma}$ is not closed.
This curve is also the concatenation of three curves belonging to $A_0\cup B_0$ which is a finite set.
 Let 
\[ C_{hor}:=\max_{\gamma} ||g(\gamma)|| \]
be the maximum word length of the corresponding group elements $g(\gamma)$, where $\gamma$ varies over this finite collection 
of loops just constructed. 
We conclude that 
\[|\Psi(\wtU_1)-\Psi(\wtU_2)|_\Sigma \leq C_{hor}.\]

The case of vertically adjacent vertices proceeds similarly in spirit.  
Suppose\footnote{Here, it is convenient to index things using zero and one, instead of one and two as above.}  
$\wtU_1, \wtU_0$ are vertically adjacent in $\Gamma$, with $|\wtU_1| = |\wtU_0|+1=n+1$, and let $e=(\wtU_1, \wtU_0)$. 
The curves $\alpha_{\wtU_1}$
and $\alpha_{\wtU_2}$ join their centers to elements $\pi_{\Sigma}(w)$ and $\pi_{\Sigma}(v)$, $w,v\in V(\Sigma)$.
Set \[\wtgamma =  \lambda_w * \alpha_{\wtU_1} * \eta_e * \overline{\alpha}_{\wtU_0}\]
which is a curve joining to points from $\pi_{\Sigma}(V(\Sigma))$ at the same level $n$.
It follows that $\gamma= f^n(\wtgamma)$ is a loop based at $b$ of the form
\[ \gamma = \lambda_x * \alpha_1 * \eta * \overline{\alpha}_0.\]
Hence it belongs to a finite set of curves (see Figure 1). 

\begin{figure}[h]
\label{fig:part3}
\includegraphics{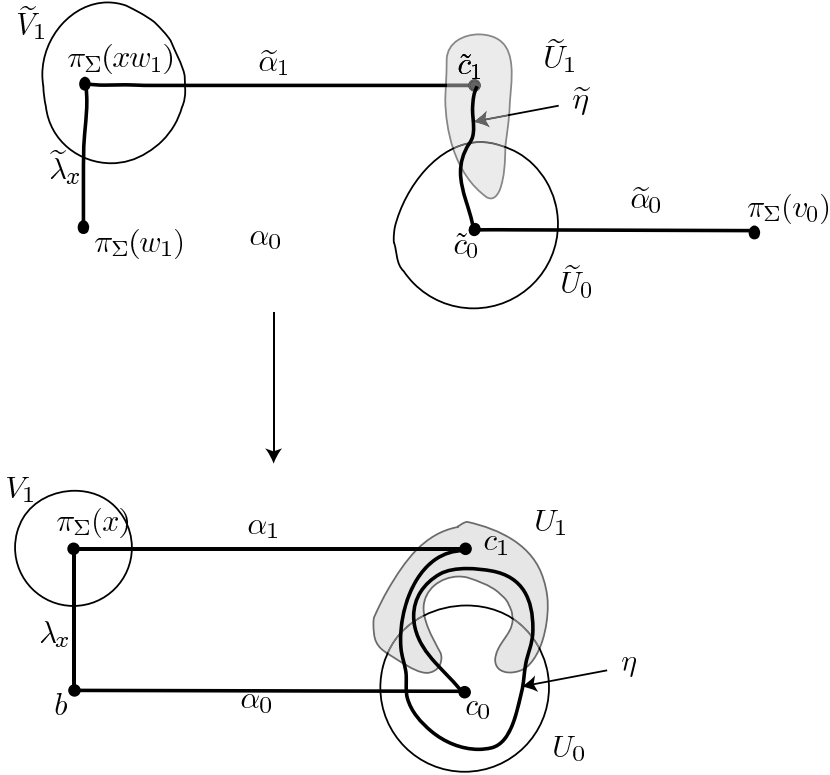}
\caption{Showing $|\Psi(\wtU_1)-\Psi(\wtU_0)|_\Sigma = O(1)$.  The disk $f(U_1)$ is not shown.}
\end{figure}

We let 
\[ C_{ver}:=\max_{\gamma} ||g(\gamma)|| \]
be the maximum word length of the corresponding group elements $g(\gamma)$, where $\gamma$ varies over the finite collection 
of loops just constructed. Thus 
\[  |\Psi(\wtU_1) - \Psi(\wtU_0)|_\Sigma \leq C_{ver}.\]

We conclude that $\Psi$ is $\max\{C_{hor}, C_{ver}\}$-Lipschitz, and the proof is complete. 
\qed

Suppose now that $f$ has no periodic critical points, and that the virtual endomorphism $\phi$ is contracting. 
Then $\Sigma$ is hyperbolic; Theorem \ref{thm:coincidence} implies that $\Gamma$ is also hyperbolic. Hence if $\epsilon$ is chosen sufficiently small, the Floyd boundary $\bdry\Gamma$ coincides with the Gromov boundary, and $\Phi$ induces a homeomorphism $\Phi: \JJJ \to \bdry \Gamma$ which conjugates the shift $\sfs: \JJJ \to \JJJ$ to the induced map $f_\Gamma: \bdry \Gamma \to \bdry \Gamma$. It follows that 
\[ \pi:= \Phi^{-1} \circ \pi_\Gamma: S^2 \to \JJJ\]
gives a semiconjugacy of branched coverings from $f: S^2 \to S^2$ to $\sfs: \JJJ \to \JJJ$.  The remainder of the proof will use finiteness properties of the dynamical system $\sfs: \JJJ \to \JJJ$  to analyze the possibilities for the fibers of $\pi$.

We conclude this section with a technical result needed for the analysis in the next section. 

By \cite[Thm. 6.13]{kmp:ph:cxcii}, $\sfs: \JJJ \to \JJJ$ is a branched covering of degree $d$.   Given $E \subset S^2$ a continuum and a component $\wtE$ of $f^{-1}(E)$, we define 
\[ \deg(f:\wtE \to E) = \min\{ \deg(f: \wtU \to U): E \subset U\}\]
where $U$ is a connected open neighborhood of $E$.  

\begin{lemma}
\label{lemma:degree}
Suppose $f_\Sigma(\txi) = \xi$, and set $E = \pi^{-1}(\xi), \wtE = \pi^{-1}(\txi)$.  Then $\deg(f:\wtE \to E) = \deg(f_\Sigma, \txi)$.  In particular, the inverse image under $f$ of a fiber of $\pi$ is a disjoint union of fibers of $\pi$, and the degrees of $f$ on these fibers sum to $d$.  
\end{lemma}

\pf We first prove $\deg(f:\wtE \to E) \geq \deg(f_\Sigma, \txi)$.  Let $U$ be a small open connected neighborhood of $E$ and $\wtU$ the pullback containing $\wtE$.  Since the decomposition of $X$ by the fibers of $\pi$ is upper-semicontinuous, and $f_\Sigma$ is a branched covering, we may assume that the images $V=\pi(U), \wtV=\pi(\wtU)$ are neighborhoods of $\xi, \txi$, respectively, that $U=\pi^{-1}(V), \wtU=\pi^{-1}(\wtV)$, and that $\deg(f_\Sigma,\txi) = \deg(f_\Sigma: \wtV \to V)$. 
Since $f_\Sigma: \JJJ \to \JJJ$ is a branched covering, the set of branch values is  nowhere dense, so there exists $\zeta \in V$ which is not a branch value of $f_\Sigma$.  Hence $\#f_\Sigma^{-1}(\zeta)\intersect \wtV = \deg(f_\Sigma,\txi)$.   It follows easily that $\deg(f_\Sigma, \txi)$ is a lower bound  for $\deg(f: \wtU \to U)$.  

The proof of the inequality $\deg(f:\wtE \to E) \geq \deg(f_\Sigma, \txi)$ is similar.  Suppose $V, \wtV$ are respectively neighborhoods of $\xi, \txi$ and set $U = \pi^{-1}(V), \wtU = \pi^{-1}(\wtV)$.  Choose $\zeta \in V$ which is not a branch value of $f_\Sigma$.  Then $\#f_\Sigma^{-1}(\zeta)=d$, and since $f$ is a branched covering of degree $d$, this implies that for each $z\in \pi^{-1}(\zeta)$, $\#f^{-1}(z)=d$ also.  Hence no two elements of $f^{-1}(z)$ belong to the same fiber of $\pi$.  We conclude $\#f_\Sigma^{-1}(\zeta)\intersect \wtV = \#f^{-1}(z)\intersect \wtU$ and so $\deg(f: \wtU \to U)$ is a lower bound for $\deg(f_\Sigma, \txi)$. 

\qed

\section{Proof of necessity}


\begin{prop} Let $f:S^2\to S^2$ be a postcritically finite topological cxc map of the sphere.
Then its virtual endomorphism is contracting.\end{prop}

Let $\OOO$ be the orbifold modeled on $(S^2, P)$, $\nu: \OOO \to \N$ its weight function, and as usual choose a basepoint $b\in S^2\setminus P$. 
Taking an iterate if necessary, we may assume
that $f(b)=b$. The group $G=\pi_1(\OOO,b)$
acts naturally as the deck transformation group of a Galois (regular) ramified universal orbifold covering $p:\widetilde{\OOO}\to S^2$ such that 
$\deg(p, \tilde{x}) = \nu(p(\tilde{x}))$ for all $\tilde{x}\in\widetilde{\OOO}$.

It follows that there exists a lift of ``$f^{-1}$'', that is, a map 
${\widetilde{f}_{-}}:\widetilde{\OOO}\to\widetilde{\OOO}$ such that
$f\circ p \circ{\widetilde{f}_{-}}= p$.  
Since $f(b)=b$, we may choose ${\widetilde{f}_{-}}$ so that a given $\tilde{b} \in p^{-1}(b)$ is fixed by ${\widetilde{f}_{-}}$: 
 \[ 
\begin{array}{ccc}
(\widetilde{\OOO}, \tilde{b}) & \stackrel{{\widetilde{f}_{-}}}{\longleftarrow} &(\widetilde{\OOO}, \tilde{b})
\\
p \downarrow & \; & \downarrow p \\
(\OOO, b)& \stackrel{f}{\longrightarrow} &(\OOO, b) \\
\end{array}
\] 
To see this: at lower-left, replace $\OOO$ with $\OOO_1$, so that the bottom arrow becomes an orbifold covering. 
Lift the orbifold structure on $\OOO_1$ at lower-left under $p$ to obtain a new orbifold structure $\widetilde{\OOO}_1$ 
on the underlying space of $\widetilde{\OOO}$ so that $p: \widetilde{\OOO}_1 \to \OOO_1$ is an orbifold covering. 
The map $\widetilde{f}_{-}$ is then a universal orbifold covering of $\widetilde{\OOO}_1$.

In this paragraph, we define a metric $\tilde{d}$ on $\OOO$ with respect to which ${\widetilde{f}_{-}}$ becomes a uniform contraction.  Let $\UUU_0$ be a covering as in \S\ref{scn:covering}, and 
endow $S^2$ with the metric $d_{\ep}$ described in \S\,\ref{scn:covering}.  
Choose $M>0$ so small that (i) any $d_{\ep}$-ball on $S^2$ of radius $2M$ is contained 
in some element of $\UUU_0$; from the choice of $\UUU_0$ it follows that 
(ii) any such ball meets $P$ in at most one point. 
Let $\DDD$ denote the set of Jordan domains $U$ of diameter at most $M$. 
Condition (ii) implies that if $U \in \DDD$, then the inverse image $p^{-1}(U)$ 
is a collection $\{\wtU_k\}$ of pairwise disjoint Jordan domains in $\widetilde{\OOO}$ 
invariant under the covering group action of $G$ 
(their closures might meet, but we don't care), and the restriction of $p$ to any domain 
$\wtU_k$ is proper. 
Let $\tilde{U}$ be one such component and set, for $x,y\in \tilde{U}$,
$$q_{\tilde{U}}(x,y)= \inf \diam_\ep p(K)$$ 
where the infimum is over all continua $K\subset \overline{\tilde{U}}$ joining $x$ to $y$. 
Now suppose $x,y$ in $\widetilde{\OOO}$.  We define the metric $\tilde{d}$ on $\widetilde{\OOO}$ by 
$$\tilde{d}(x,y)=\inf \sum_{j} q_{\tilde{U}_j} (x_j,x_{j+1})$$
where the infimum is taken over every chain $x=x_1, x_2, \ldots, x_\ell = y$ joining $x$ to $y$, with the property that for each $j$, there exists $U_j \in \DDD$ and a component $\wtU_j$ of $p^{-1}(U_j)$   such that $\{x_j, x_{j+1}\}\subset \wtU_j$.  
We note that if $\tilde{d}(x,y)\le M$, then we may find $\tilde{U}\in\DDD$ such that $\tilde{d}(x,y)=q_{\tilde{U}}(x,y)$.
It follows that the infimum in the definition of $\tilde{d}$ is actually attained. 

The function $\tilde{d}$ defines a metric and the group $G$ acts by isometries on $(\widetilde{\OOO},\tilde{d})$.

\begin{lemma} 
\label{lemma:orbqi}
The map $\Phi:g\mapsto g(\tilde{b} )$ defines a quasi-isometry from $G$ to $(\widetilde{\OOO},\tilde{d})$.\end{lemma}

\pf We already know that $G$ is finitely generated and that it
acts on $\widetilde{\OOO}$ by isometries, properly discontinuously and cocompactly. 

We mimic the proof of the Svarc-Milnor lemma. 
Let $C\subset \widetilde{\OOO}$ be a compact set containing $\tilde{b}$ such that $G.C=\widetilde{\OOO}$
and fix $D\ge 8M$ such that $C\subset B(\tilde{b} , D/4)$ and let $S=\{g\in G, \tilde{d}(g(\tilde{b} ),\tilde{b})\le D\}$.

Then $\Phi(G)$ is cobounded, $S$ is a finite generating set and there are constants $\lambda$, $c>0$ 
such that $\tilde{d}(\Phi(g_1),\Phi(g_2))\le \lambda \|g_1g_2^{-1}\|_S + c$
for all $g_1,g_2\in G$, see \cite[Prop.\,8.19]{bridson:haefliger:book}. The only
point to prove is the reverse inequality.

Fix $g\in G$ and let us consider $x_0,\ldots, x_n$ with $x_0=\tilde{b} $ and $x_n=g(\tilde{b} )$ and 
let $\tilde{U_1},\ldots, \tilde{U}_n$
be such that $$\tilde{d}(\tilde{b} ,g(\tilde{b} ))=\sum q_{\tilde{U}_j} (x_{j-1},x_j)\,.$$

Let $N$ be the integer part of $4\tilde{d}(\tilde{b} ,g(\tilde{b} ))/D$. For any $1\le k \le N$, let $y_k$ be the first $x_i$ such that
$$\sum_{j=1}^i q_{\tilde{U}_j} (x_{j-1},x_j) \ge k D/4\,.$$ It follows that 
$\tilde{d}(y_k,y_{k+1})\le D/2$ since the points $(x_j)$ are at most $M$-separated and $M\le D/8$.
Choose now $g_k\in G$ such that $\tilde{d}(g_k(\tilde{b} ),y_k)\le D/4$; then $\tilde{d}(g_k(\tilde{b} ),g_{k+1}(\tilde{b} ))\le D$, so that
$g_k\circ g_{k+1}^{-1}\in S$ and $\|g\|_S\le N+1$.  Therefore, $\|g\|_S\le (4/ D ) (\tilde{d}(\tilde{b} ,g(\tilde{b} )) +1)+1$.
The proof is complete.\qed

\gap

\pf (Proposition) 
Let $x,y\in\widetilde{\OOO}$ be such that $\tilde{d}(x,y)\le M$, then we may find a continuum $K$ containing $\{x,y\}$
such that $\tilde{d}(x,y)=\diam_\ep p(K)$. Note that $p(K)$ is contained in an element of $\UUU_0$.
Therefore, for any $n\ge 0$, $(p \circ {\widetilde{f}_{-}^{n}})(K)$ is contained in an element of $\UUU_n$ and
so $\tilde{d}( {\widetilde{f}_{-}^n}(x),{\widetilde{f}_{-}^n}(y))\lesssim e^{-\ep n}$.
It follows easily that, for any pair of points $x,y$ without any restriction on their distance, we have
$$\tilde{d}( {\widetilde{f}_{-}^n}(x),{\widetilde{f}_{-}^n}(y))\lesssim e^{-\ep n} \tilde{d}(x,y)\,.$$

This implies that
$$\limsup_{n\to\infty}\left( \sup_{x,y\in\widetilde{\OOO}}
\frac{ \tilde{d}( {\widetilde{f}_{-}^n}(x),{\widetilde{f}_{-}^n}(y)) }{\tilde{d}(x,y)}\right)^{1/n} = e^{-\ep} < 1.$$

Now,  since $f(b)=b$ is fixed, we may define the virtual endomorphism $\phi$ by choosing the connecting path $\al$ 
(as in the definition given in \S 3.1) to be the constant path.
Hence, $g\in  \mbox{\rm dom}  \phi^n$ if and only if ${\widetilde{f}_{-}^n}(g(\tilde{b} ))$ belongs to the orbit  $G.\tilde{b} $,  
and in this case, $\Phi(\phi^n(g))= {\widetilde{f}_{-}^n}(g(\tilde{b} ))$. 
From the preceding Lemma \ref{lemma:orbqi}, the metrics $|| \cdot ||$ and $\tilde{d}$ are comparable, so we obtain
$$\limsup_{n\to\infty}\left( \limsup_{||g||\to\infty}
\frac{||\phi^n(g)||}{||g||}\right)^{1/n}  \le \limsup_{n\to\infty}\left( \sup_{x,y\in\widetilde{\OOO}}
\frac{ \tilde{d}( {\widetilde{f}_{-}^n}(x),{\widetilde{f}_{-}^n}(y)) }{\tilde{d}(x,y)}\right)^{1/n} < 1.$$

\qed


\section{Proof of sufficiency}

From the introduction, recall that the key step is to establish 
\gap

\begin{prop}\label{prop:fibers} Let $f:S^2\to S^2$ be a Thurston map with postcritical set $P$ and with hyperbolic biset $\frak{M}$.
The semiconjugacy $\pi: S^2 \to \JJJ$ is injective on $P$, and no fiber of $\pi$ separates $S^2$.  
\end{prop}

The main ingredients in the proof of the claim are \cite[Lemma 6.11, Theorem 6.15]{kmp:ph:cxcii}, which assert that the dynamical system $f_\Sigma: \JJJ \to \JJJ$ is expanding, irreducible, and of {\em topologically finite type}.   Finite type implies that there is a positive integer $p$ such that for all $\xi \in \JJJ$ and all $n > 0$, $\deg(f_\Sigma^n, \xi) \leq p$.   Thus in particular $f_\Sigma$ has no periodic branch points.  Irreducible in this context is equivalent to the fact that given any nonempty open set $W$, there is an integer $n$ with $f_\Sigma^n(W)=\JJJ$.

We remark that these properties imply that $f_\Sigma: \JJJ \to \JJJ$, as a topological dynamical system, is a so-called finite type topologically coarse expanding conformal (cxc) system; cf. \cite{kmp:ph:cxci}. 
This forces strong constraints on the fibers of $\pi$.  We prove the claim by a case-by-case analysis of the possibilities. 
We remark that there is much overlap with the flavor of the arguments in \cite{kmp:tan:rmwdjs}.  

\gap

A continuum $E \subset S^2$ we call 
\bi
\item {\em special} if $E \intersect P\neq \emptyset$;
\item {\em type III} if $S^2 - E$ has at least $ 3$  components that intersect $P$
\item {\em type II} if $S^2 - E$ has exactly $2$  components that intersect $P$
\item {\em type I} if $S^2 - E$ has exactly $1$  component that intersects $P$
\ib

The map $f$ sends fibers of $\pi$ onto fibers of $\pi$. Suppose $E=\pi^{-1}(\xi)$ for some $\xi\in \JJJ$.
\bi
\item[(i)]  If $E$ is special, then $f^k(E)$ is special for all $k\ge 0$, and since $P$ is finite, $E$ is eventually periodic.
Moreover, there are only finitely many such fibers.
\item[(ii)]   It is easily verified that if $f(E)$ is not special, then neither is $E$, and the  type of $f(E)$ is at least the type of $E$.   
\item[(iii)]  If $E$ is type II  and non-special then it determines a homotopy class of unoriented simple closed curve in $S^2-P$.  
\item[(iv)]  Since distinct fibers are disjoint, there can be only finitely many type III fibers, and only finitely many homotopy classes of curves arise from type II fibers.   
\item[(v)] If $U$ is a component of $S^2\setminus E$, then $\pi(U)$ is a component of $\JJJ\setminus\{\xi\}$ and is open. This follows
from the facts that $\pi$ is monotone, continuous and onto  and that  $S^2$ is locally connected.
\ib

The proof of the following technical lemma is straightforward.

\begin{lemma}
\label{lemma:thickening}
Given a continuum  $E \subset S^2$, there is a ``thickening'' $\widehat{E}$  of $E$ given by an open regular neighborhood of $E$ with the following properties.
\bi
\item the closure of $\widehat{E}$ is a compact surface with boundary. 
\item $\widehat{E}$ has the same type as $E$.
\item $\widehat{E} \intersect P = E \intersect P$.
\item for each component $\widetilde{E}$ of $f^{-1}(E)$, we have $\deg(f: \widetilde{\widehat{E}} \to \widehat{E}) = \deg(f: \widetilde{E} \to E)$, 
where $\widetilde{\widehat{E}}$ is the component of $f^{-1}(\widehat{E})$ containing $\widetilde{E}$.
\ib
\end{lemma}

In what follows, given a fiber $E=\pi^{-1}(\xi)$, the symbol $\widehat{E}$ denotes a thickening of $E$ satisfying the conditions of Lemma \ref{lemma:thickening}.  The following paragraphs rule out certain types of fibers. 
\gap

\begin{lemma}[no periodic separating fibers]\label{lma:periodicfibers}
No periodic $E$ separates $S^2$, and there are no type II fiber homotopic rel $P$ to one of its iterates. Moreover, if
$E$ is periodic and special, then $\# E \intersect P= 1$ and $E\cap P$ is periodic of the same period
as $E$.
\end{lemma}

\pf  Taking an iterate if necessary, we may assume that $E$ is fixed. 
Since $f_\Sigma$ has no periodic branch points, $\deg(f_\Gamma, \xi)=\deg(f: E \to E)=1$, by Lemma \ref{lemma:degree}. 
Then $E$ is contained in a thickening $\widehat{E}$ for which $\deg(\widetilde{\widehat{E}} \to \widehat{E})=1$. 
This implies $f|_E$ is a homeomorphism which extends to a homeomorphism on a neighborhood of $E$. 

\be

\item {\bf Subcase $E$ special, $\# E \intersect P\geq 2$.}  Since $f:E\to E$ is a homeomorphism and $P$ is forward invariant and finite, we may assume that
$P\cap E$ is pointwise fixed. Note that $\hat{E}\cap f(\hat{E})$ is a neighborhood of $E$. 
Let $\gamma\subset f(\hat{E})$ join two distinct point of $E\cap P$. Then $f^{-1}\circ\gamma$ joins the same points 
and is homotopic rel. $P$ to a curve $\gamma_1\subset f(\hat{E})$  
since $f|E$ is $1-1$. 
This implies the existence of an arc with arbitrarily long nontrivial iterates 
under the pullback relation, which is impossible by Theorem \ref{thm:prevents_levy}. 

\item {\bf Subcase $E$  special and type II or type III.} Let us first establish that 
the boundary $C$ of some complementary component $W$ of $\widehat{E}$ must be essential and nonperipheral in $S^2-P$.  

If  each such $W$ contains at most one point of $P$ then $f^{-1}(\widehat{E})$ is connected. Thus
$\deg(f: f^{-1}(\widehat{E}) \to \widehat{E})=d>1$ holds, and this contradicts the fact that $\deg(f: E \to E)=1$.  
Hence at least one $W$ contains two points of $P$. If $E$ is type III or type II and special, then $C$ is essential and nonperipheral.

We first note that $f^k:C\to f^k(C)$ has degree $1$ up to homotopy rel. $P$ for all $k\ge 1$. Moreover,  $f^k(C)$ is also essential and nonperipheral
for otherwise, $C$ would not be. Finally, the forward orbit of $C$ has to cycle eventually up to homotopy since there are only finitely homotopic
classes of fibers of type at least II.  This
means that a Levy cycle can be extracted from $\{f^k(C)\}_{k\ge 0}$, contradicting Theorem \ref{thm:prevents_levy}.

\item {\bf Subcase $E$  non special and type II or type I.} Let us assume that $E$ separates $S^2$. 
There ia a component $W$ of $S^2\setminus E$
which contains all of $P$ but at most one single point.   Let 
$\widehat{W}$ be the component of $S^2\setminus \widehat{E}$  in $W$.  Then the complement 
$\widehat{U}=S^2\setminus \widehat{W}$ is a Jordan domain containing $E$ and at most one critical value of $f$.
Therefore, $f^{-1}(\widehat{U})$ is a union of Jordan domains.
One of them  contains $E$, hence  $U= S^2\setminus W$. It follows that $f(U)\subset U$. But $U$ contains a component
of $S^2\setminus E$, so $\pi(U)$ has interior by (v) and $\sfs(\pi(U))\subset \pi(U)$. which contradicts the irreducibility of $\sfs$.

\item {\bf Subcase type II $E$ periodic up to homotopy.} We assume that $E$ is type II and $f^k(E)$ is homotopic
to $E$ rel. $P$. If $\deg(f^k, E)>1$ then $\deg(f^{nk},E) =\deg(f_\Sigma^{nk}, \xi) \to \infty$ as $n \to \infty$, violating $f_\Sigma$ being finite type.  
So $f^k|_E$ is injective, and after thickening, we extract a Levy cycle as above. 
\eb\qed

\gap

{\noindent\bf No preperiodic separating fibers nor type II and type III fibers.} 
Let $E$ be a preperiodic separating fiber. If all its iterates are nonspecial, then we would obtain
a periodic separating fiber by (ii) and (iv), which is impossible according to Lemma \ref{lma:periodicfibers}.
Therefore, there is some $n\ge 1$ such that $E'=f^n(E)$ is special and periodic by (i).
By Lemma \ref{lma:periodicfibers} again, $E'$ does not separate $S^2$ and intersects $P$ at a single point, so
the inclusion $U:=S^2-(P\union E')\hookrightarrow S^2-P$ is a homotopy equivalence.  
Therefore $f^{-n}(U)$ is homotopy equivalent to the connected set 
$S^2-f^{-1}(P)$.  Hence $U$ is connected, $E$ is nonseparating, and $\#E \intersect P \leq 1$.  

Let us now prove that any  type II and type III fibers have to cycle up to homotopy, which will establish their nonexistence.. 
Since fibers cannot intersect, there can be only finitely many of them up to homotopy, cf. (iv).
Let $E$ be of type at least II. If its iterates are non special, then they remain of type at least II by (ii), and hence have to cycle up to homotopy by (iv). 
Otherwise, one iterate becomes special, hence is preperiodic as well. 

\gap

In order to complete the proof, we must prove that there are no type I fiber. We have already ruled out preperiodic ones.
\gap

\noindent{\bf No wandering separating type I fiber.}  All the possible remaining separating type I fibers are nonspecial and hence form a totally invariant set of $S^2$.
We use the same notation as in Lemma \ref{lma:periodicfibers}: given such a fiber $E$, let $\widehat{W}_E\subset W_E$
contain  $P$, $U_E=S^2\setminus W_E$ and
$\widehat{U}=S^2\setminus \widehat{W}$. 
Therefore, $f^{-1}(\widehat{U}_{f(E)})$ is a union of Jordan domains disjoint from $P$.
One of them  contains $E$, hence  $U_E$. It follows that $f(U_E)\subset U_{f(E)}$. 
This implies that no preimage of $P$ can enter the sets $U_E$: this contradicts again the irreducibility of $\sfs$. 
Therefore, no type I fiber separates the sphere. 

\gap

\pf (Proposition \ref{prop:fibers}) It follows from the discussion above that 
fibers of $\pi$ are connected, do not separate the sphere, and meet $P$ in at most one point.   \qed

\pf (Theorem \ref{thm:characterization})  Suppose the virtual endomorphism $\phi$ is contracting and $f$ has no periodic branch points.  
Necessarily  $\#P \geq 3$.   For technical reasons, it is convenient later to assume $\#P > 3$ (see the introduction for the case $\#P=3$). 

According to Theorem \ref{thm:coincidence}, we may consider defining data for $\Sigma$ and $\Gamma$ so that both graphs are quasi-isometric.
Since $\phi$ is contracting, it follows that $\Gamma$ is hyperbolic and that the quasi-isometry defines a homeomorphism
$\Phi:\JJJ\to\partial\Gamma$. The map $\pi=\Phi^{-1}\circ\pi_\Gamma:S^2\to \JJJ$ semiconjugates $f$ to $\sfs$. 
Moreover, Lemma \ref{lemma:monotone} and Proposition \ref{prop:fibers} imply that  $\pi$  is
monotone, injective on $P$ and that its fibers do not separate $S^2$.
By Moore's Theorem \cite[Theorem 25.1]{daverman:decompositions}, the quotient $\JJJ$ is homeomorphic to the sphere. 
The induced map $f_\Sigma: \JJJ \to \JJJ$ is thus an expanding Thurston map.

It remains to construct a homotopy from $f$ to an expanding map $g$. By \cite[Theorems 13.4, 25.1]{daverman:decompositions}, 
the decomposition $\GGG$ of $S^2$ by the fibers of the semiconjugacy $\pi: S^2 \to \JJJ$ has the property of being {\em strongly shrinkable}: 
there is a one-parameter family of continuous maps $h_t: S^2 \to S^2, t \in [0,1]$, called a {\em pseudoisotopy}, such that $h_0 = \id_{S^2}$, 
$h_t$ is a homeomorphism for $t \in [0,1)$, $h_t|_{P_f} = \id_{P_f}$ for all $t$, and the fibers of $h_1$ and those of $\pi$ coincide.  
The induced homeomorphism $S^2 \to \JJJ$ conjugates $\sfs$ to an expanding Thurston map $g: S^2 \to S^2$, and 
the family of maps $f_t$, $t \in [0,1]$ defined as the unique continuous solution of $h_t \circ f = f_t \circ h_t$, with $f_1=g$, gives a  homotopy
through Thurston maps. 
The proof of the sufficiency in 
Theorem \ref{thm:characterization} is complete. \qed


\end{document}